\documentclass[]{article}
\usepackage{amssymb,amsmath,amsthm,mathrsfs,bm,setspace}

 \newcommand{\R}{\mathbf{R}}

\renewcommand{\P}{\mathrm{P}} \newcommand {\E}{\mathrm{E}}
\newcommand{\1}{{\bf 1}} \renewcommand{\d}{\text{\rm d}}
\newcommand{\e}{\text{\rm e}}

\newcommand{\bdim}{\text{\rm B-}\underline{\dim}}
\newcommand{\Bdim}{\text{\rm B-}\overline{\dim}}

\newtheorem{stat}{Statement}[section]
\newtheorem{proposition}[stat]{Proposition}

\newtheorem{theorem}[stat]{Theorem}  \newtheorem{lemma}[stat]{Lemma}
\theoremstyle{definition}
\newtheorem{definition}[stat]{Definition}\newtheorem{remark}[stat]{Remark}

 \numberwithin{equation}{section}
\title{\bf On dynamical bit sequences}
\date{June 8, 2009}
\author{%
  Davar Khoshnevisan\\University of Utah \and David
  A. Levin\\University of Oregon \and Pedro
  J. M\'endez-Hern\'andez\\Universidad de Costa Rica}

\begin{document}
\onehalfspacing
\maketitle
\begin{abstract}
  Let $\bm{X}_k(t):=(X_1(t),\ldots,X_k(t))$ denote a $k$-vector of
  i.i.d.\ random variables, each taking the values 1 or 0 with
  respective probabilities $p\in(0\,,1)$ and $1-p$. As a process
  indexed by $t\ge 0$, $\bm{X}_k$ is constructed---following
  Benjamini, H{\"a}ggstr{\"o}m, Peres, and Steif \cite{BHPS}---so that
  it is strong Markov with invariant measure
  $((1-p)\delta_0+p\delta_1)^k$. We derive sharp estimates for the
  probability that ``$X_1(t)+\cdots+X_k(t)=k-\ell$ for some $t\in
  F$,'' where $F\subseteq[0\,,1]$ is nonrandom and compact. We do this
  in two very different settings: (i) Where $\ell$ is a constant; and
  (ii) Where $\ell=k/2$, $k$ is even, and $p=q=1/2$. We prove that the
  probability is described by the Kolmogorov capacitance of $F$ for
  case (i) and Howroyd's $\frac12$-dimensional box-dimension profiles
  for case (ii). We also present sample-path consequences, and a
  connection to capacities that answers a question of \cite{BHPS}.
	
	
\end{abstract}

\section{Introduction and main results}
      
Choose and fix some $p\in(0\,,1)$. By a \emph{bit sequence} we mean
simply a $k$-vector $(X_1,\ldots,X_k)$---where $k\ge 1$ is fixed---of
independent, identically-distributed random variables with
$\P\{X_1=1\}=1-\P\{X_1=0\}=p$. For simplicity, we write $\P_p$ in
order to keep track of all dependencies on the parameter $p$.

A \emph{dynamical bit sequence} is the process $\bm{X}_k
:=\{(X_1(t)\,,\ldots,X_k(t)\}_{t\ge 0}$ that is constructed as follows: 
We begin with a bit sequence $(X_1\,,\ldots,X_k)$ at time zero. Then,
to every index $j\in\{1\,,\ldots,k\}$ we associate a rate-one Poisson process;
all Poisson processes being independent of all $X$ variables.
And whenever the $j$th Poisson process jumps, we replace the
corresponding variable $X_j$ by a copy, independent of all else,
in order to obtain a time-dependent family of random variables.
More precisely,  let  $\{X_i^{(j)}\}_{1 \leq i,j < \infty}$ be an i.i.d.\ array of random bits,
each satisfying $\P_p\{X_i^{(j)} = 1\} = p$, and consider also
a sequence of independent rate-one Poisson processes, independent
of the $\{X_i^{(j)}\}_{1 \leq i,j < \infty}$. 
Denote the jump times of the $i$-th Poisson process by
$J^{(0)}_i := 0 < J^{(1)}_i < J^{(2)}_i < \cdots,$
and define
\begin{equation}
	X_i(t) := \sum_{j=1}^\infty X_i^{(j)}\1_{\left\{
	J_i^{(j)}\le t < J_i^{(j-1)}\right\}}
	\quad\text{for $t\ge 0$ and $1\le i\le k$}.
\end{equation}

Dynamical bit sequences were introduced by Benjamini, H{\"a}ggstr{\"o}m, Peres,
and Steif in 2003. They can be used to model how bit sequences can get
corrupted over time, for instance, see \cite{BHPS} and  
\cite{H:DPS,HPS}. A closely-related variant of dynamical bit sequences
was introduced earlier by Rusakov \cite{Rusakov,Rusakov:94,Rusakov:93}; see also 
Rusakov and Chuprunov \cite{RC}.

A key feature of $\bm{X}_k$ is that it is a strong Markov process on
$\{0\,,1\}^k$ whose invariant measure is the law
$(p\delta_1+(1-p)\delta_0)^k$ of the bit sequence
$(X_1\,,\ldots,X_k)$; in particular, $(X_1(t)\,,\ldots,X_k(t))$ and
$(X_1\,,\ldots,X_k)$ have the same distribution for every fixed $t\ge
0$.

The main goal of the present work is to estimate a family of
probabilities for dynamical bit sequences.  Among other
things, these estimates can be used to describe almost-sure properties
of runs [i.e., contiguous sequences of ones and/or zeros]; see
\cite{BHPS} and Section \ref{sec:runs} below for examples.

Let $S_k:=X_1+\cdots+X_k$.  We note the classical binomial
identity
\begin{equation}\label{eq:bin1}
  \P_p\left\{ S_k=k-\ell\right\}
  =\binom{k}{\ell}p^{k-\ell}(1-p)^\ell,
\end{equation}
valid for all integers $\ell=0,\dots,k$. Consequently,
\begin{equation}\label{eq:bin2}
  \P_p\left\{ S_k=k-\ell\right\}
  \asymp k^\ell p^k\quad\text{as $k\to\infty$},
\end{equation}
for every fixed integer $\ell\ge 0$, where ``$a_k\asymp b_k$ for a
given set of values of $k$'s'' means that ``$(a_k/b_k)$ is bounded
above and below by positive and finite constants, uniformly for all
mentioned values of $k$.''

In contrast to \eqref{eq:bin2}, one can verify that for every fixed
integer $\ell\in\{0\,,\ldots,k\}$,
\begin{equation}\label{eq:BHPS}
  \P_p\left\{ S_k(t)=k-\ell\text{ for some $t\in[0\,,1]$}\right\}
  \asymp k^{\ell+1} p^\ell,
\end{equation}
as $k\to\infty$, where $S_k(t):=X_1(t)+\cdots+X_k(t)$.  Benjamini et
al \cite{BHPS} have proved \eqref{eq:BHPS} in the case that $\ell=0$,
and the more general case follows from their methods as
well. Moreover, the time set $[0\,,1]$ can be replaced by any other
time interval without affecting the form of \eqref{eq:BHPS}.

Our main goal is to describe the effect of the geometry of the time
variable $t$ on the bounds in \eqref{eq:BHPS}.  In order to describe
our first main estimate, choose and fix a compact set
$F\subseteq\R_+$. Then for all $\epsilon>0$ define ${\rm
  K}_F(\epsilon)$ to be the largest integer $m\ge 1$ for which there
exist points $t_1,\ldots,t_m\in F$ such that $|t_i-t_j|\ge\epsilon$
for all $i\neq j$. The function ${\rm K}_F$ is known as the
\emph{Kolmogorov $\epsilon$-capacity} [or \emph{capacitance}, or
\emph{$\epsilon$-packing}] of the set $F$. Now we can describe or
first result.

\begin{theorem}\label{th:main1}
  Choose and fix an integer $\ell\ge 1$ and a nonrandom compact set
  $F\subset\R_+$. Then as $k\to\infty$,
  \begin{equation}\label{eq:main1}
    \P_p\left\{\exists t\in F:\ S_k(t)
      =k-\ell\right\}
    \asymp {\rm K}_F(1/k)\,k^\ell p^k.
  \end{equation}
\end{theorem}

When $F=\{0\}$ this recovers \eqref{eq:bin2} because ${\rm
  K}_{\{0\}}(1/k)=1$.  And when $F=[0\,,1]$, it recovers
\eqref{eq:BHPS} since $k{\rm K}_{[0,1]}(1/k)\to 1$ as
$k\to\infty$. The following implies more interesting behavior when $F$
is less regular.

\begin{remark}
  Let
  \begin{equation}
    \alpha:=\liminf_{k\to\infty}
    \frac{\log {\rm K}_F(1/k)}{\log k}
    \quad\text{and}\quad
    \beta:=
    \limsup_{k\to\infty}\frac{\log {\rm K}_F(1/k)}{\log k}	\end{equation}
  respectively denote the lower,
  and the upper, Minkowski dimension of $F$. Then
  it follows from Theorem \ref{th:main1} that
  \begin{equation}
    k^{\ell+\alpha+o(1)}p^k\le
    \P_p\left\{ \exists t\in F:
      S_k(t)=k-\ell\right\}
    \le k^{\ell+\beta+o(1)} p^k.
  \end{equation}
  And each bound is achieved along a suitable subsequence of $k$'s.
  \qed
\end{remark}

The probability in Theorem \ref{th:main1} has even more interesting
behavior in the cases that $\ell\to\infty$ at a prescribed rate. We
investigate one critical case in this paper.  First, let us apply
Stirling's formula in \eqref{eq:bin1} to recover the well-known local
limit theorem $\P_{1/2} \{ S_k= k/2 \}\asymp k^{-1/2}$, valid for all
sufficiently-large even integers $k$.  Our next result describes the
significantly-different dynamical version of this estimate, and is
motivated by sequence-matching estimates of Arratia and Waterman
\cite{AW:89,AW:85} in the what turns out to be the critical case where
$p=1/2$ and $\ell=k/2$.  Let
\begin{equation}
  \gamma:=\bdim_{1/2} F
  \quad\text{and}\quad
  \delta:=\Bdim_{1/2}F
\end{equation}
respectively denote Howroyd's $\frac12$-dimensional lower, and upper,
box-dimension profile of $F$ \cite{Howroyd:01}. [These dimension
numbers will be recalled later on during the course of the proof of
Theorem \ref{th:main2} below; see \S\ref{sec:3}.]  Then we have the
following.

\begin{theorem}\label{th:main2}
  If $k\to\infty$ along the even integers, then for every nonrandom
  bounded set $F\subset\R_+$,
  \begin{equation}\label{eq:main2}
    k^{-(1-\gamma)/2+o(1)}\le
    \P_{1/2}\left\{ \exists t\in F:\ S_k(t) =k/2\right\}
    \le k^{-(1-\delta)/2+o(1)}.
  \end{equation}
  Moreover, each bound is achieved along a suitable subsequence of
  $k$'s.
\end{theorem}
Thus, in the cases where $\gamma=\bdim_{1/2}F<\Bdim_{1/2}F=\delta$, we
can deduce that the probability in \eqref{eq:main2} decays roughly
like a power of $k$, and nevertheless there are no power-law
asymptotics.

Theorems \ref{th:main1} and \ref{th:main2} are proved respectively in
\S\ref{sec:T1} and \S\ref{sec:3}. We collect some consequences of
Theorem \ref{th:main1} in \S\ref{sec:runs}.  Finally, we describe a
connection to Riesz capacity in \S\ref{sec:cap}, and use it to verify
a conjecture of Benjamini et al \cite{BHPS}.

\section{Proof of Theorem \ref{th:main1}}\label{sec:T1}

\begin{lemma}\label{lem:Z}
  For every integers $\ell=0\,,\ldots,k$ and real number $t\ge 0$,
  \begin{equation}\begin{split}
      &\P_p\left( S_k(t)=k-\ell \ \big|\, S_k(0)=k-\ell\right)\\
      &\hskip1in=\sum_{i=0}^{\min(\ell,k-\ell)} \binom{\ell}{i}
      \theta_t^i (1-\theta_t)^{\ell-i}\binom{k-\ell}{i}
      \kappa_t^{k-\ell-i} (1-\kappa_t)^i,
    \end{split}\end{equation}
  where $\theta_t := p(1-\e^{-t})$ and $\kappa_t := p + \e^{-t}q$.
\end{lemma}

\begin{proof}
  
  We first find the transition probabilities
  for the process  $\{X_1(t)\}_{t \geq 0}$.  Let $\{P_t\}_{t\ge 0}$
  define 
  the transition matrices defined by
  \begin{equation}
    P_t(i\,,j) 
      := \P_p\{X_1(t) = j \mid X_1(0) = i\} \quad
           \text{for all $i,j \in \{0\,,1\}$} \,.  
  \end{equation}

  Recall that  $\{J_1^{(k)}\}_{k=1}^\infty$ are
  the jump times of the rate-one Poisson  clock associated  to $X_1$. 
  Then   
  \begin{equation} \begin{split} \label{Eq:P1}
     \P_p\bigl\{\,X_1(t) = 1 \mid & X_1(0)  = 1\,\bigr\} \\
       & = \frac{1}{p} \; \P_p \left\{\, X_1(t) = 1, X_1(0) = 1, 
          t< J_1^1\,\right\} \\
       & \quad  + \frac{1}{p} \; \P_p \left\{\, X_1(t) = 1
		, X_1(0) = 1, t \geq J_1^1 \, \right\} \\
       & = \frac{1}{p} \;\left[\, p e^{-t} + p^2 (1-e^{-t}) \right],
  \end{split} \end{equation}
  
  and this quantity is manifestly equal to $\kappa_t$.

  In fact, we can follow the same argument to conclude  that 
  \begin{equation} \label{eq:Pt}
    P_t =
    \begin{bmatrix}
      pe^{-t} + q & p(1-e^{-t}) \\
     q(1-e^{-t}) & qe^{-t} + p \\
    \end{bmatrix}
    =
    \begin{bmatrix}
      1-\theta_t & \theta_t \\
      1- \kappa_t & \kappa_t
    \end{bmatrix}
    \quad\text{for all $t\ge 0$.}
  \end{equation}

	
  For all $u,v\in\{0\,,1\}$ define
  \begin{equation}
    N_{u\to v} := \sum_{j=1}^{k-1}
    \1_{\{X_j(t)=v,X_j(0)=u\}}.
  \end{equation}
  This quantity denotes the number of integers
  $j\in\{1\,,\ldots,k-1\}$ such that $X_j(0)=u$ and yet $X_j(t)=v$. It
  follows from the strong Markov property and \eqref{eq:Pt} that the
  following properties are valid under the conditional measure
  $\P_p(\,\cdots\,|\, S_k(0)=k-\ell)$:
  \begin{itemize}
  \item $N_{0\to 1}$ has the binomial distribution with parameters
    $\ell$ and $\theta_t$;
  \item $N_{1\to 1}$ has the binomial distribution with parameters
    $k-\ell$ and $\kappa_t$;
  \item $N_{0\to 1}$ and $N_{1\to 1}$ are independent.
  \end{itemize}
  Additionally, the conditional probability in the statement of the
  lemma is
  \begin{equation}\begin{split}
      &\P_p\left\{ N_{0\to 1}+N_{1\to 1}=k-\ell
        \,\big|\, S_k(0)=k-\ell\right\}\\
      &\hskip.6in= \sum_{i=0}^{\min(\ell,k-\ell)}
      \P_p\left\{ N_{0\to 1}=i\,\big|\, S_k(0)=k-\ell\right\}\\
      &\hskip1.8in \times \P_p\left\{ N_{1\to 1}=k-\ell-i\,\big|\,
        S_k(0)=k-\ell\right\}.
    \end{split}\end{equation}
  The lemma follows from these observations.
\end{proof}

The proof of Theorem \ref{th:main1} uses a 2-scale argument that is
borrowed from our earlier work on dynamical random walks \cite{KLM};
it can be outlined as follows: First we prove that if $I$ is a
``small'' closed interval, then the two events $\{S_k(t)=k-\ell\
\text{for some $t\in I$}\}$ and $\{ S_k(0)=k-\ell\}$ have more or less
the same chances of occurring (Proposition \ref{pr:hit}). We will also
demonstrate that ``small'' means ``whose length is of sharp order
$1/k$''; this length scale---or correlation length---was found earlier
in \cite{BHPS}.  Then we cover our set $F$ with closed intervals of
length $1/k$, and apply a covering argument. Finally, we show that
this covering argument produces sharp answers. With this outline in
mind, we begin with our first step.

\begin{proposition}\label{pr:hit}
  As $k\to\infty$,
  \begin{equation}
    \P_p\left\{ S_k(t)=k-\ell
      \text{ for some $t\in[0\,,1/k]$}\right\}
    \asymp k^\ell p^k.
  \end{equation}
\end{proposition}

\begin{proof}
  Let $\mathcal{F}:=\{\mathcal{F}_t\}_{t\ge 0}$ denote the filtration
  generated by the strong Markov process $\bm{X}_k$.  We can assume,
  without loss of generality, that $\mathcal{F}$ is augmented in the
  usual manner. Define
  \begin{equation}
    L_k := \int_0^{2/k} \1_{\{S_k(t)=k-\ell\}}\,\d t.
  \end{equation}
  According to \eqref{eq:bin2},
  \begin{equation}\label{eq:EI}
    \E_p(L_k) \asymp k^{\ell-1} p^k
    \quad\text{as $k\to\infty$}.
  \end{equation}
  Next, we consider the stopping times, $\sigma:= \inf\{
  t\in[0\,,1/k]:\ S_k(t)=k-\ell\}$, where
  $\inf\varnothing:=\infty$. By the strong Markov property, the
  following holds almost sure $[\P_p]$ on $\{\sigma<\infty\}$:
  \begin{equation}\begin{split}
      \E_p\left( L_k\,\big|\, \mathcal{F}_\sigma\right) &=
      \E_p\left(\left. \int_\sigma^{2/k} \1_{\{S_k(t)=k-\ell\}}
          \, \d t\ \right|\, \mathcal{F}_\sigma\right)\\
      &\ge \int_0^{1/k} \P_p\left\{\left.
          S_k(t+\sigma)=k-\ell\ \right|\, \mathcal{F}_\sigma\right\}\,\d t\\
      &= \int_0^{1/k} \P_p\left\{ S_k(t)=k-\ell\, \big|\,
        S_k(0)=k-\ell\right\}\,\d t.
    \end{split}\end{equation}
  We apply Lemma \ref{lem:Z} to find that $\P_p$-almost surely on
  $\{\sigma<\infty\}$,
  \begin{equation}
    \E_p\left( L_k\,\big|\, \mathcal{F}_\sigma\right)
    \ge \int_0^{1/k} (1-\theta_t)^a \kappa_t^{k-\ell}\,\d t
    \ge\frac 1k\int_0^1\left( 1-\theta_u\right)^k \kappa_{u/k}^k
    \,\d u.
  \end{equation}
  The integrand converges to $\e^{-u}$ as $k\to\infty$. Consequently,
  thanks to the bounded convergence theorem,
  \begin{equation}
    \E_p\left( L_k\,\big|\, \mathcal{F}_\sigma\right)
    \ge \frac{\text{const}}{k}\cdot\1_{\{\sigma<\infty\}}
    \qquad\text{$\P_p$-almost surely}.
  \end{equation}
  Here, the implied constant does not depend on either $k$ or $\ell$.
  We can take expectations of both side, solve, and apply the optional
  stopping theorem of Doob to find that $\P_p\{\sigma<\infty\}\le
  \text{const}\cdot k\E_p (L_k)$.  The proposition follows from this
  and \eqref{eq:EI}.
\end{proof}

\subsection*{Proof of Theorem \ref{th:main1}}
We divide the proof in to two parts. The first part (the upper bound)
is a relatively simple covering argument. In the second part, we show
that the covering argument produces the correct answer.  For that we
use a second-moment argument.

\begin{proof}[Proof of the upper bound in \eqref{eq:main1}]
  Let $\mathcal{J}_k$ denote the collection of all closed subintervals
  $I$ of $[0\,,1]$ that have the form $I=[i/k\,,(i+1)/k]$, where
  $i\in\{0\,,\ldots,k-1\}$. According to Proposition \ref{pr:hit},
  \begin{equation}\begin{split}
      \P_p\left\{\exists t\in F:\, S_k(t)=k-\ell\right\}&\le
      \sum_{\substack{I\in\mathcal{J}_k\\ I\cap F\neq\varnothing}}
      \P_p\left\{\exists t\in I:\, S_k(t)=k-\ell\right\}\\
      &\le \text{const}\cdot\left|\left\{ I\in\mathcal{J}_k:\, I\cap
          F\neq\varnothing\right\}\right|\cdot k^\ell p^k,
    \end{split}\end{equation}
  where $|\cdots|$ denotes cardinality. The lemma follows since the
  preceding cardinality is $\le 3{\rm K}_F(1/k)$ \cite[Proposition
  2.7]{KLM}.  The upper bound in \eqref{eq:main1} follows readily from
  this.
\end{proof}

\begin{proof}[Proof of the lower bound in \eqref{eq:main1}]
  According to the definition of ${\rm K}_F(1/k)$, we can find
  $ t_1<\cdots<t_{{\rm K}_F(1/k)} $ in $F$ such that $ t_{i+1}-t_i\ge
  k^{-1} $ for  $ 1 \leq i \leq {\rm K}_F(1/k)-1 $.  It follows easily that
  \begin{equation}\label{eq:KolmPoint}
    t_j-t_i\ge \frac{j-i}{k}\qquad
    \text{for $1\le i<j\le{\rm K}_F(1/k)$}.
  \end{equation}
	
  Define
  \begin{equation}
    N_k := \sum_{j=1}^{{\rm K}_F(1/k)} \1_{\{S_k(t_j)=k-\ell\}}.
  \end{equation}
  According to \eqref{eq:bin2},
  \begin{equation}\label{eq:EJ}
    \E_p(N_k) \asymp {\rm K}_F(1/k)\, k^\ell p^\ell
    \quad\text{as $k\to\infty$}.
  \end{equation}
  Observe that $\E_p(N_k^2)$ is equal to
  \begin{equation}\label{eq:N2}
    \E_p(N_k) + \mathop{\sum\sum}\limits_{%
      1\le i\neq j\le {\rm K}_F(1/k)} \P_p\left\{
      S_k(|t_j-t_i|)=k-\ell\,, S_k(0)=k-\ell \right\}.
  \end{equation}

  Define
  \begin{equation}\label{eq:T}
    \mathcal{T}:=\P_p\left\{ S_k(h)=k-\ell\ \big|\,
      S_k(0)=k-\ell\right\}.
  \end{equation}
  By the Markov property and Lemma \ref{lem:Z},
  \begin{equation}\label{eq:PS}
    \mathcal{T}=\sum_{i=0}^{\min(\ell,k-\ell)}Q^{(0)}_i(h),
  \end{equation}
  where for all $h> 0$ and $i\in\{0\,,\ldots,\min(\ell\,,k-\ell)\}$,
  \begin{equation}\begin{split}
      Q^{(0)}_i(h) &:= \binom{\ell}{i}\binom{k-\ell}{i}\theta_h^i
      (1-\theta_h)^{\ell-i}\kappa_h^{k-\ell-i}
      (1-\kappa_h)^i\\
      &\le \ell^{(\ell)} k^i \theta_h^i(1-\theta_h)^{\ell-i}
      \kappa_h^{k-\ell-i} (1-\kappa_h)^i.
    \end{split}\end{equation}
  In order to simplify this object, let us first define
  \begin{equation}
    a:=\max(p\,,1-p)
    \quad\text{and}\quad
    b:=\min(p\,,1-p).
  \end{equation}
  Let us also observe that
  \begin{equation}
    \theta_h(1-\kappa_h)=ab\left(1-\e^{-h}\right)^2
    \le \frac{a}{b}\left[ 1-\left(a+b\e^{-h}\right)\right]^2.
  \end{equation}
  And
  \begin{equation}
    (1-\theta_h)^{\ell-i}\kappa_h^{k-\ell-i} \le
    \left(a+b\e^{-h}\right)^{k-2i}.
  \end{equation}
  It follows that uniformly for all $\ell\in\{0\,,\ldots,k\}$ and
  $i\in\{0\,,\ldots,\min(\ell\,,k-\ell)\}$,
  \begin{equation}\label{eq:Q0Q1}
    Q^{(0)}_i(h) \le \text{const}\cdot k^i Q^{(1)}_i(h),
  \end{equation}
  where
  \begin{equation}
    Q^{(1)}_i(h) := \left[ 1-\left(a+b\e^{-h}\right)\right]^{2i}
    \left(a+b\e^{-h}\right)^{k-2i}.
  \end{equation}
  The function $[0\,,1]\ni\lambda\mapsto
  (1-\lambda)^{k-2i}\lambda^{2i}$ is maximized at $\lambda:=2i/k$, and
  the value of the maximum is at most $k^{-2i}$. And therefore, by
  \eqref{eq:Q0Q1},
  \begin{equation}\label{eq:Q1}
    Q^{(0)}_i(h) \le \text{const}\cdot k^{-i}. 
  \end{equation}
  The preceding is a good estimate when $i\ge 1$. In the case that
  $i=0$, we merely use
  \begin{equation}
    Q^{(0)}_0(h) =(1-\theta_h)^\ell\kappa_h^{k-\ell}
    \le (a+b\e^{-h})^k.
  \end{equation}
  We plug this and \eqref{eq:Q1} in \eqref{eq:PS}, using \eqref{eq:T},
  to find
  \begin{equation}\label{eq:PS1}\begin{split}
      \mathcal{T}&\le\text{const}\cdot \left(
        (a+b\e^{-h})^k+\sum_{i=1}^{\min(\ell,k-\ell)} k^{-i}\right)\\
      &\le\text{const}\cdot \left( (a+b\e^{-h})^k+ k^{-1}\right).
    \end{split}\end{equation}
  Therefore, it follows from \eqref{eq:bin2}, \eqref{eq:KolmPoint},
  and \eqref{eq:N2} that $\E_p(N_k^2)$ is at most
  \begin{align}
      &\E_p(N_k) + \text{const}\cdot k^\ell p^k
      \mathop{\sum\sum}\limits_{%
        1\le i\neq j\le {\rm K}_F(1/k)} \left(
        \left(a+b\e^{-|t_j-t_i|}\right)^k+ k^{-1}\right)\\\nonumber
      &\qquad\le \E_p(N_k) + \text{const}\cdot k^\ell p^k
      \mathop{\sum\sum}\limits_{%
        1\le i\neq j\le {\rm K}_F(1/k)} \left(
        \left(a+b\e^{-|j-i|/k}\right)^k+ k^{-1}\right).
    \end{align}
  Consequently, $\E_p(N_k^2)$ is bounded from above by
  \begin{align}
	&\E_p(N_k)+\text{const}\cdot {\rm K}_F(1/k)\cdot k^\ell
		p^k\sum_{j=1}^{{\rm K}_F(1/k)}\left(
		\left(a+b\e^{-j/k}\right)^k+ k^{-1}\right) \\\nonumber
	& \le\E_p(N_k)+\text{const}\cdot {\rm K}_F(1/k)
		\cdot k^\ell p^k  \left\{ k \int_0^1 (a+b\e^{-u})^k\,\d u +\frac{{\rm
		K}_F(1/k)}{k}\right\}.
	\end{align}
  [The final term is at most one because ${\rm K}_F(1/k) \le {\rm
    K}_{[0,1]}(1/k)=k$.] A direct computation of the preceding
  integral reveals that
  \begin{equation}
    \E_p(N_k^2)\le \E_p(N_k) + \text{const}\cdot {\rm K}_F(1/k)\,k^\ell p^k
    \le \text{const}\cdot \E_p(N_k);
  \end{equation}
  see \eqref{eq:EJ}. Therefore, \eqref{eq:EJ} and the Paley--Zygmund
  inequality \cite[p.\ 72]{Kh} together imply that $\P_p\{N_k>0\} \ge
  \text{const}\cdot {\rm K}_F(1/k)\,k^\ell p^k.$ The theorem follows
  because $\{N_k>0\} \subseteq\{\exists t\in F:\, S_k(t)=k-\ell\}$.
\end{proof}

\section{Proof of Theorem \ref{th:main2}}\label{sec:3}

Let $\mathcal{P}(F):=$ all Borel probability measures $\mu$ such that
$\mu(F)=1$.

\begin{theorem}\label{th:ell=a/2}
  If $k$ is an even integer, then for every measurable and nonrandom
  set $F\subseteq[0\,,1]$, the following holds as $k\to\infty$:
  \begin{equation}\begin{split}
      &\P_{1/2} \left\{ \exists t\in F:\
        S_k(t)=k/2\right\} \\
      &\hskip.7in\asymp \frac{1}{\sqrt k}\cdot \left[
        \inf_{\mu\in\mathcal{P}(F)}\iint\min\left(
          \frac{1}{\sqrt{k|t-s|}}~,~1\right) \,\mu(\d s)\,\mu(\d
        t)\right]^{-1}.
    \end{split}\end{equation}
\end{theorem}

As a first key step, we simplify the expression in Lemma \ref{lem:Z}.

\begin{lemma}\label{lem:cond:a/2}
	Uniformly for every even integer $k\ge 1$ and $t\in[0\,,1]$:
	\begin{equation}
		\P_{1/2}\left( S_k(t)=k/2 \ \big|\ S_k(0)=k/2\right)
		\asymp \min\left(
		\frac{1}{\sqrt{kt}}~,~ 1\right).
	\end{equation}
\end{lemma}

\begin{proof}
  Throughout, we will use the following simple fact: 
  \begin{equation}\label{eq:exp}
    \frac{t}{2} \leq  1-\e^{-2t} \le 2t
    \qquad\text{for all $t\in[0\,,1]$}.
  \end{equation}

 We appeal to the result, and notation, of Lemma 2.1.  Notice that in
  the present case, $\theta_t=1-\kappa_t=(1-\e^{-t})/2$, and moreover,
  \begin{equation}
    \P_{1/2}\left(  S_k(t)=k/2 \ \big|\ S_k(0)=k/2\right)
    =\P\{ X=Y\},
  \end{equation}
  where $X$ and $Y$ are two independent binomial random variables with
  common parameters $k/2$ and $(1-\e^{-t})/2$.  We can apply the
  Plancherel formula:
  \begin{equation}
    \P\{X=Y\} =\frac{1}{2\pi}\int_{-\pi}^\pi
    |\phi(z)|^2\, \d z,
  \end{equation}
  where $\phi$ is the characteristic function of $X$; i.e.,
  \begin{equation}
    \phi(z) = \left( \frac{1-\e^{-t}}{2}+\frac{1+\e^{-t}}{2} \e^{i z}\right)^{k/2}.
  \end{equation}
  Of course,
  \begin{equation}
    \left|\frac{1-\e^{-t}}{2}+\frac{1+\e^{-t}}{2} \e^{i z}\right|^2=
    \frac{1+\e^{-2t}}{2}+\frac{1-\e^{-2t}}{2}\cos z.
  \end{equation}

 Taylor's theorem with remainder shows that $1-(z^2/2)\le \cos z \le
  1-(z^2/20)$ for all $z\in[-\pi\,,\pi]$.  Therefore,
  \begin{equation}
    1-\frac{z^2}{4}\left(1-
      \e^{-2t}\right)\le
    \left|\frac{1-\e^{-t}}{2}+\frac{1+\e^{-t}}{2} \e^{i z}\right|^2 
    \le  1 - \frac{z^2}{40}\left(1-\e^{-2t}\right).
  \end{equation}
  Now,  (\ref{eq:exp}) implies that for all $t\in[0\,,1]$
  and $z\in[-\pi\,,\pi]$,
  \begin{equation}
    1-\frac{z^2t}{2}\le
    \left|\frac{1-\e^{-t}}{2}+\frac{1+\e^{-t}}{2} \e^{i z}\right|^2 
    \le  1 - \frac{z^2t}{80};
  \end{equation}
  and therefore, 
  \begin{equation}
    \int_{-1}^1 \left(\,1-\frac{z^2t}{2}\,\right)^{k/2}\,\d z
    \le \int_{-\pi}^\pi|\phi(z)|^2\,\d z \le 
    \int_{-\pi}^\pi \,\left(\,1-\frac{z^2t}{80}\,\right)^{k/2}\d z.
  \end{equation}
  Since $ 1-u\leq \e^{-u}$ for all $u\ge 0$, we find that 
  uniformly for all $t\in[0\,,1]$,
  \begin{equation}
     \P\{X=Y\}  \le 
    \frac{1}{2\pi} \int_{-\pi}^\pi \e^{-z^2tk/160}\,\d z 
    \le\frac{\text{const}}{\sqrt{tk}}.
  \end{equation}
  This implies the upper bound of the lemma. For the lower bound follows 
  from the fact $ 1-u/2\ge \e^{-u}$, for all $0\le u \le 1$, and
  $\int_{-1}^1\e^{-z^2tk/2}\,\d z\ge\text{const}\cdot
  \min((tk)^{-1/2}\,,1)$, uniformly for all $t\in[0\,,1]$.
\end{proof}

Next we prove Theorem \ref{th:ell=a/2}.

\begin{proof}[Proof of Theorem \ref{th:ell=a/2}]
  For all $\mu\in\mathcal{P}(F)$ and even integers $k\ge 1$ define
  \begin{equation}
    L_k^\mu := \int \mathbf{1}_{\{S_k(t)=k/2\}}\, \mu(\d t).
  \end{equation}
  In accord with \eqref{eq:EI}, the following is valid: As
  $k\to\infty$,
  \begin{equation}\label{eq:EL}
    \E_{1/2} (L_k^\mu) \asymp k^{-1/2},
  \end{equation}
  where the approximation holds uniformly for all $\mu\in\mathcal{P}
  (F)$. Next we estimate the second moment of $L_k^\mu$, using Lemma
  \ref{lem:cond:a/2}:
  \begin{equation}\begin{split}
      \E_{1/2}\left[ \left( L_k^\mu\right)^2\right] &\le 2
      \mathop{\iint}\limits_{t\ge s} \P_{1/2}\left\{
        S_k(t)=S_k(s)=k/2\right\}\,
      \mu(\d s)\,\mu(\d t)\\
      &\le\frac{\text{const}}{\sqrt k}\cdot\iint \min\left(
        \frac{1}{\sqrt{k|t-s|}} ~,~ 1\right)\,\mu(\d s)\,\mu(\d t).
    \end{split}\end{equation}
  This, and the Paley--Zygmund inequality, together imply that
  \begin{equation}
    \P_{1/2} \{ L_k^\mu>0 \}\ge \frac{\text{const}}{\sqrt k}\cdot
    \left[ \iint
      \min\left( \frac{1}{\sqrt{k|t-s|}} ~,~ 1\right)\,\mu(\d s)\,
      \mu(\d t)\right]^{-1}.
  \end{equation}
  The preceding probability is at most $\P_{1/2}\{\exists t\in F:\,
  S_k(t)=k/2\}$. Since the latter does not depend on
  $\mu\in\mathcal{P}(F)$,
  \begin{equation}\begin{split}
      &\P_{1/2}\left\{\exists t\in F:\, S_k(t)=k/2\right\}\\
      &\hskip.3in\ge\frac{\text{const}}{\sqrt k}\cdot\left[\inf_{
          \mu\in\mathcal{P}(F)} \iint\min\left(
          \frac{1}{\sqrt{k|t-s|}} ~,~ 1\right)\,\mu(\d s)\, \mu(\d
        t)\right]^{-1}.
    \end{split}\end{equation}
  This proves half of the theorem.
	
  For a converse bound, let us consider the stopping time
  $\sigma:=\inf\{t\in F:\, S_k(t)=k/2\}$, where
  $\inf\varnothing:=2$. We apply our existing computation with the
  following special choice of $\mu$: $\mu(\bullet) := \P_{1/2} (
  \sigma\in \bullet\,|\, \sigma<1)$.  By the strong Markov property
  and Lemma \ref{lem:cond:a/2},
  \begin{equation}\begin{split}
      \E_{1/2}\left( \left. L_k^\mu\ \right|\,
        \mathcal{F}_\sigma\right) &=\int_{s\ge\sigma} \P_{1/2}\left(
        S_k(s)=k/2\,\big|\, \mathcal{F}_\sigma
      \right)\,\mu(\d s)\\
      &\ge \text{const}\cdot\int_{s\ge\sigma} \min\left(
        \frac{1}{\sqrt{k(s-\sigma)}}~,~1\right) \,\mu(\d s),
    \end{split}\end{equation}
  where the implied constant does not depend on $\mu$, and is
  nonrandom.  Since $\sigma$ is a bounded stopping time, we can apply
  the optional stopping theorem to deduce from the preceding that for
  all probability measures $\mu$ on $F$,
  \begin{equation}
    \E_{1/2}(L_k^\mu)\ge\text{const}\cdot
    \E_{1/2} \left[\int_{s\ge\sigma} \min\left(
        \frac{1}{\sqrt{k(s-\sigma)}}~,~1\right)
      \,\mu(\d s)\right].
  \end{equation}
  According to \eqref{eq:EL},
  \begin{equation}
    \E_{1/2} \left[\int_{s\ge\sigma} \min\left(
        \frac{1}{\sqrt{k(s-\sigma)}}~,~1\right)
      \,\mu(\d s)\right] \le \frac{\text{const}}{\sqrt k}.
  \end{equation}
  On the other hand, the definition of $\mu$ implies that the
  left-hand side is equal to
  \begin{equation}\begin{split}
      &\mathop{\iint}\limits_{s\ge t}\min\left(
        \frac{1}{\sqrt{k(s-t)}}~,~1\right)
      \,\mu(\d s)\,\mu(\d t)\cdot \P_{1/2}\{\sigma<1\}\\
      &\hskip.4in\ge\frac12\iint\min\left(
        \frac{1}{\sqrt{k|t-s|}}~,~1\right) \,\mu(\d s)\,\mu(\d t)\cdot
      \P_{1/2}\{\sigma<1\}.
    \end{split}\end{equation}
  It follows from the preceding two displays that
  \begin{equation}
    \P_{1/2}\{\sigma<1\} \le \frac{\text{const}}{\sqrt k}\cdot
    \left[\iint\min\left(
        \frac{1}{\sqrt{k|t-s|}}~,~1\right)
      \,\mu(\d s)\,\mu(\d t)\right]^{-1}.
  \end{equation}
  This proves the theorem.
\end{proof}

Finally, we prove Theorem \ref{th:main2}.

\begin{proof}[Proof of Theorem \ref{th:main2}]
  We recall Howroyd's theory of box-dimension profiles
  \cite{Howroyd:01}: For every $s>0$ and $x\in\R$ define
  $\psi_s(x):=\min(1\,,|x|^{-s})$.  Then, given $r>0$, a sequence of
  pairs $(w_i\,,x_i)_{i=1}^n$ is a \emph{size-$r$ weighted
    $\psi_s$-packing of $F$} if: (i) $x_i\in F$ for all $i$; (ii)
  $w_i\ge 0$ for all $i$; and (iii) $\sum_{j=1}^n w_j\psi_s
  ((x_i-x_j)/r)\le 1$, uniformly for all $i=1,\ldots,n$.  Define
  $N_r(F\,;\psi_s) := \sup\sum_{i=1}^n w_i$, where the supremum is
  taken over all size-$r$ weighted $\psi_s$-packings
  $(w_i\,,x_i)_{i=1}^n$ of $F$.  The \emph{$s$-dimensional upper
    box-dimension profile} $\Bdim_s F$ of $F$ is then defined as
  \begin{equation}
    \Bdim_s F := \limsup_{r\downarrow 0}\frac{\log N_r(F\,;\psi_s)}{
      \log(1/r)},
  \end{equation}
  where $\log 0:=-\infty$.  The \emph{$s$-dimensional lower
    box-dimension profile} $\bdim_s F$ of $F$ is defined as above, but
  with a $\liminf$ in place of the $\limsup$.
	
  According to Khoshnevisan and Xiao \cite[Theorem 4.1]{KX},
  \begin{equation}
    \Bdim_s F=\limsup_{n\to\infty}
    \frac{1}{\log (1/n)}\log \inf_{\mu\in\mathcal{P}(F)}
    \iint \psi_s\left(\frac{u-v}{1/n}\right)\,\mu(\d u)\,\mu(\d v).
  \end{equation}
  That proof shows also that $\bdim_s F$ is equal to the same
  quantity, but with $\limsup$ replaced by $\liminf$.  In light of
  these two facts, Theorem \ref{th:main2} follows from Theorem
  \ref{th:ell=a/2}.
\end{proof}

\section{Some applications to runs}\label{sec:runs}

Throughout this section we will be studying $\bm{X}_\infty$ and its
dynamical version; both can be defined in the usual way via infinite
product spaces, as respective ``limits'' of $\bm{X}_k$ and its
dynamical version, as $k\to\infty$. We skip the details, as they are
outlined nicely in \cite{BHPS}.

As was pointed out in \cite{BHPS}, one can use an estimate such as
\eqref{eq:BHPS} to study the behavior of the longest runs in a
dynamical bit sequence, as the sequence length tends to infinity. We
describe this work next. For every two integers $n,\ell\ge 1$, define
$Z_n^{(\ell)}$ to be the largest integer $j\ge\ell+1$ such that
$X_m=1$ for all but $\ell$ values of $m\in\{n\,,\ldots,n+j-1\}$. If
such a $j$ does not exist, then $Z_n^{(\ell)}:=0$.  One defines the
dynamical version similarly: $Z_n^{(\ell)}(t)$ is defined as above, by
$X_m$ is replaced by $X_m(t)$.

According to the Erd\H{o}s--R\'enyi theorem \cite{ErdosRenyi}, for
every $\ell\ge 1$,
\begin{equation}
  \lim_{n\to\infty} \frac{Z_n^{(\ell)}}{\log_{1/p}n}=1
  \qquad\text{$\P_p$-almost surely},
\end{equation}
where $\log_{1/p}$ denotes the base-$(1/p)$ logarithm. Erd\H{o}s and
R\'ev\'esz \cite{ErdosRevesz} improved this statement by showing that
the following holds for every nonrandom sequence
$\bm{a}:=\{a_j\}_{j=1}^\infty$ of positive integers that tend to
infinity:
\begin{equation}\label{eq:ER}
  \P_p\left\{ \limsup_{n\to\infty} Z_n^{(\ell)}\ge a_n\text{ i.o.}
  \right\}=\begin{cases}
    0&\text{if $\sum_{n=1}^\infty a_n^\ell p^{a_n}<\infty$},\\
    1&\text{otherwise}.
  \end{cases}
\end{equation}
This particular formulation appears explicitly, for example, in the
book by R\'ev\'esz \cite[p.\ 60]{Revesz} in the case that $p=1/2$.  To
be more precise, R\'ev\'esz (\emph{loc.\ cit.})  defines
$Z_n^{(\ell)}$ as the longest run having \emph{at most} $\ell$ defects
in the first $n$ bits. But a real-variable comparison argument reveals
that our definition and that of R\'ev\'esz have the same asymptotic
behavior \cite{ErdosRevesz}.

Of course, $Z_n^{(\ell)}$ can be replaced by $Z_n^{(\ell)}(t)$ for a
fixed $t$, and \eqref{eq:ER} continues to holds. By contrast, it is
shown in \cite[Theorem 1.4]{BHPS} that
\begin{equation}\label{eq:R}
  \P_p\left\{\exists t\in [0\,,1]:\ Z_n^{(\ell)}(t)\ge a_n
    \text{ i.o.} \right\}=\begin{cases}
    0&\text{if $\sum_{n=1}^\infty a_n^{\ell+1} p^{a_n}<\infty$},\\
    1&\text{otherwise}.
  \end{cases}
\end{equation}
To be precise, Benjamini et al \cite{BHPS} prove this for $\ell=0$;
but one can apply their method, using \eqref{eq:BHPS}, to produce
\eqref{eq:R}. We learned about \eqref{eq:R} from the monograph by
R\'ev\'esz \cite[p.\ 61]{Revesz}, who conjectured \eqref{eq:R} in the
case that $p=1/2$.  It is possible to study the size of the set of
times $t$ at which $Z^{(\ell)}_n(t)\ge a_n$ infinitely often.  Define
\begin{equation}\label{eq:EDefn}
  \mathcal{E}^{(\ell)}(\bm{a}) := \left\{ t\in[0\,,1]:\
    Z_n^{(\ell)}(t)\ge a_n\text{ i.o.}\right\}.
\end{equation}
In light of \eqref{eq:R}, $\mathcal{E}^{(\ell)}(\bm{a})$ is nonempty
if and only if $\sum_{n=1}^\infty a_n^{\ell+1}p^{a_n}=\infty$.  Then,
it is possible to adapt the argument of \cite[Theorem 1.5]{BHPS},
using \eqref{eq:BHPS} of the present paper, to derive the following:
\begin{equation}
  \dim_{_{\rm H}}\mathcal{E}^{(\ell)}(\bm{a})=\sup\left\{
    s\in(0\,,1):\ \sum_{n=1}^\infty a_n^{\ell+1-s}p^{a_n}
    =\infty\right\},
\end{equation}
where $\dim_{_{\rm H}}$ denotes Hausdorff dimension.

Equation \eqref{eq:ER} is a statement about whether or not
$\mathcal{E}^{(\ell)}(\bm{a})$ is void. Our next result describes all
nonrandom compact sets $F\subset[0\,,1]$ that have positive
probability of intersecting $\mathcal{E}^{(\ell)}(\bm{a})$.  First we
need some notation from \cite{KLM}.

We say that an interval is \emph{rational} if its endpoints are
rational numbers. From here on, we choose and fix a nonrandom set $F$
and a nonrandom sequence $\bm{a} :=\{a_j\}_{j=1}^\infty$ of positive
integers that tend to infinity.  It turns out that the following
definition \cite{KLM} defines the correct intersection property for
the random set $\mathcal{E}^{(\ell)}(\bm{a})$:
\begin{definition}
  We write ``$\Psi(\bm{a}\,;F)<\infty$'' if and only if we can find
  closed rational intervals $F_1,F_2,\ldots$ such that: (i)
  $\cup_{n=1}^\infty F_n\supseteq F$; and (ii)
  \begin{equation}
    \sum_{j=1}^\infty {\rm K}_{F_n}(1/a_j)
    a_j^\ell p^{a_j}<\infty
    \quad\text{for all $n\ge 1$.}
  \end{equation}
\end{definition}

\begin{theorem}\label{th:hit:a.s.}
  For every nonrandom compact set $F\subseteq[0\,,1]$,
  \begin{equation}
    \P_p\left\{ \mathcal{E}^{(\ell)}(\bm{a})\cap F\neq\varnothing\right\}
    =\begin{cases}
      0&\text{if $\Psi(\bm{a}\,;F)<\infty$},\\
      1&\text{otherwise}.
    \end{cases}
  \end{equation}
\end{theorem}

One can adapt the proof of \cite[Theorem 2.5]{KLM} to the present
setting to see that the following implies Theorem \ref{th:hit:a.s.}.
We will not prove Theorem \ref{th:hit:a.s.}, since its proof does not
require new ideas. However, we will prove the following crucial step,
whose proof borrows ideas from the proof of Theorem 1.4 of
\cite{BHPS}.

\begin{proposition}\label{pr:Mdim}
  For all nonrandom compact sets $F\subseteq[0\,,1]$,
  \begin{equation}
    \P_p\left\{ \sup_{t\in F}Z^{(\ell)}_n(t)\ge a_n
      \text{ i.o.}
    \right\}=\begin{cases}
      0&\text{if $\sum_{n=1}^\infty {\rm K}_F(1/a_n)
        a_n^\ell p^{a_n}<\infty$},\\
      1&\text{otherwise}.
    \end{cases}
  \end{equation}
\end{proposition}

\begin{proof}
  According to Theorem \ref{th:main1}, $\sum_{n=1}^\infty {\rm
    K}_F(1/a_n) a_n^\ell p^{a_n}<\infty$ if and only if
  \begin{equation}
    \sum_{n=1}^\infty \P\left\{\sup_{t\in F} Z^{(\ell)}_n(t)\ge
      a_n \right\}<\infty.
  \end{equation}
  Therefore, whenever $\sum_{n=1}^\infty {\rm K}_F(1/a_n) a_n^\ell
  p^{a_n}$ is finite, $\sup_{t\in F}Z^{(\ell)}_n(t)\ge a_n$ only
  occurs finitely often.  For the converse let us suppose $\sum_{n=1}^\infty
  {\rm K}_F(1/a_n) a_n^\ell p^{a_n}=\infty$, and define the events
  \begin{equation}
    F_n := \left\{
      X_{n-j}(t)=0\text{ for every $j=1,2,\ldots,\ell+1$
        and $t\in [0\,,1]$}\right\}.
  \end{equation}
  Note that $F_n$ is independent of the event
  \begin{equation}\begin{split}
      \tilde{F}_n&:=\left\{\sup_{t\in F}
        Z^{(\ell)}_n(t) \ge a_n\right\}\\
      &= \left\{ \exists t\in F:\, S_{n+a_n-1}(t)-S_n(t)=a_n-\ell
      \right\},
    \end{split}\end{equation}
  and $\P_p(F_n)\ge (q/\e)^{\ell+1}$ with $q:=1-p$.  Since
  $\sum_{n=1}^\infty {\rm K}_F(1/a_n) a_n^\ell p^{a_n}=\infty$,
  \begin{equation}
    \sum_{j=1}^\infty \P_p(G_j)=\infty
    \quad\text{where}\quad
    G_n := F_n\cap \tilde{F}_n.
  \end{equation}
  Note that for all integers $N\ge 1$,
  \begin{equation}
    \sum_{n,m=1}^N\P_p(G_n\cap G_m) 
    =\sum_{n=1}^N\P_p(G_n)
    +2\sum_{n,m=1}^N
    \P_p(G_n\cap G_m).
  \end{equation}
  Suppose $n<m$. If $m\le n+a_n-1$, then the events $G_n$ and $G_m$
  are disjoint. If $m\ge n+a_n+\ell+1$, then $G_n$ and $G_m$ are
  independent. In the remaining $O(1)$ cases we can use the elementary
  bound $\P_p(G_n\cap G_m)\le \P_p(G_n)$. Since
  $\sum_n\P_P(G_n)=\infty$, it follows that $\sum_{n,m=1}^N\P_p
  (G_n\cap G_m )$ is bounded above by a constant times
  $(\sum_{j=1}^N\P_p(G_j))^2$, where the implied constant does not
  depend on $N\ge 1$. The Borel--Cantelli lemma for dependent event
  \cite{ChungErdos} implies that infinitely-many of the $G_n$'s---and
  hence the events $\tilde{F}_n$---occur almost surely.
\end{proof}

Let us end this section with the following example: Let $\theta>0$ be
fixed, and consider the sequence $\bm{a}(\theta)$ given by
\begin{equation}
  a_n = a_n(\theta) := l_pn+\theta l_pl_pn,
\end{equation}
where $l_p x := \log_{1/p}(\max(x\,,100))$.  It is possible to check
that
\begin{equation}
  \sum_{n=1}^\infty {\rm K}_F(1/a_n(\theta))\,
  \left[a_n(\theta)\right]^\ell p^{a_n(\theta)}<\infty
  \quad\text{iff}\quad\sum_{n=100}^\infty
  \frac{{\rm K}_F(1/\log n)}{n(\log n)^{\theta-\ell}}\,\d s<\infty.
\end{equation}
The doubling property, ${\rm K}_F(2\epsilon) \asymp {\rm
  K}_F(\epsilon)$, valid for all $\epsilon>0$ sufficiently small
\cite[eq.\ (2.8)]{KLM}, implies that
\begin{equation}
  \sum_{n=1}^\infty {\rm K}_F(1/a_n(\theta))\,
  \left[a_n(\theta)\right]^\ell p^{a_n(\theta)}<\infty
  \quad\text{iff}\quad
  \int_1^\infty \frac{{\rm K}_F(1/s)}{s^{\theta-\ell}}\,\d s<\infty.
\end{equation}
According to Proposition 2.9 of \cite{KLM},
\begin{equation}
  \dim_{_{\rm P}}F+\ell+1=\inf\left\{
    \theta:\, \Psi(\bm{a}(\theta)\,;F)<\infty\right\},
\end{equation}
where $\dim_{_{\rm P}}$ denotes packing dimension. Therefore, we can
combine the preceding facts to deduce the following:
\begin{equation}
  \sup_{t\in F}\limsup_{n\to\infty}
  \frac{Z^{(\ell)}_n(t)-l_pn}{l_pl_pn}=\dim_{_{\rm P}}F+\ell+1
  \quad\text{$\P_p$-a.s.}
\end{equation}
When $F:=\{0\}$, this is due to Erd\H{o}s and R\'enyi
\cite{ErdosRenyi}; and when $F:=[0\,,1]$ it is due to Benjamini et al
\cite[Theorem 3.1]{BHPS}.

\section{A sharp capacity criterion}\label{sec:cap}

Let $m:\R_+\to\R_+$ be a strictly increasing function so that
$m(\mathbf{N})\subset\mathbf{N}$.  Benjamini et al \cite{BHPS} have
proposed the following ``bit process'' as part of their parity test
that is motivated by complexity theory: Define
\begin{equation}
  B_k(t):=\bigoplus_{j=m(k)}^{m(k+1)} X_j(t),
\end{equation}
where $t\ge 0$, $k\in\mathbf{N}_+$, and $\oplus$ denotes addition mod
2. Of course, $B_k(t)$ is either zero or one. It is proved in
\cite[Lemma 4.1]{BHPS} that
\begin{equation}\label{eq:2/m}
  \P_{1/2}\left\{\exists t\in[0\,,1]:\
    B_k(t)=0\text{ for all $k\in\mathbf{N}$}\right\}>0
  \quad\text{iff}\quad
  \sum_{k=1}^\infty \frac{2^k}{m(k)}<\infty,
\end{equation}
provided that $m$ satisfies the Hadamard gap condition,
\begin{equation}\label{eq:gap}
  \inf_{k\ge 1}\frac{m(k+1)}{m(k)}>1.
\end{equation}

Consider the random set
\begin{equation}\label{eq:T_m}
  T_m := \left\{t\in[0\,,1]:\ B_k(t)=0
    \text{ for all $k\in\mathbf{N}$}\right\}.
\end{equation}
Consider the special case that $m$ is the function $m_q(x):=[2^{x/q}]$
for some fixed $q>0$.  Lemma 4.1 of \cite{BHPS} shows that for all
nonrandom compact sets $E\subset[0\,,1]$,
\begin{equation}\label{eq:17}
  \text{Cap}_q(E)>0
  \quad\Rightarrow\quad
  \P_{1/2}\left\{ T_{m_q}\cap E\neq\varnothing\right\}>0
  \quad\Rightarrow\quad
  \mathcal{H}^q(E)>0,
\end{equation}
where $\text{Cap}_q$ denotes the $q$-dimensional Riesz capacity and
$\mathcal{H}^q$ the $q$-dimensional Hausdorff measure \cite[App.\ C \&
D]{Kh}. In their Remark 4.5, Benjamini et al (\textit{loc.\ cit.})
point out that there is a [small] gap between the conditions of
positive Hausdorff measure versus positive capacity. If we specialize
the next theorem to $m:=m_q$ then we obtain a verification of their
conjecture.

\begin{theorem}\label{th:cap}
  Suppose, in addition, that $t\mapsto 2^{-t}m(t)$ and $m$ are both
  strictly increasing. Additionally, choose and fix a nonrandom
  compact set $E\subset[0\,,1]$.  Then, $\P_{1/2}\{T_m\cap
  E\neq\varnothing\}>0$ iff there exists a probability measure $\rho$
  on $E$ such that $J(\rho)<\infty$, where
  \begin{equation}
    J(\rho) := \iint (\mathcal{L}g)(|t-s|)\,\rho(\d s)\,\rho(\d t);
  \end{equation}
  $g$ is defined by $\log_2 g(t)=m^{-1}(t)$, and
  $(\mathcal{L}g)(\lambda):=\int_0^\infty \e^{-\lambda s}\,g(\d s)$
  denotes the Laplace transform of the Stieltjes measure $\d g$.
\end{theorem}

\begin{remark}
  The monotonicity condition on $t\mapsto 2^{-t}m(t)$ implies the gap
  condition \eqref{eq:gap}. Indeed, under the monotonicity condition
  we have $m(k+1)/m(k)\ge 2$ for all $k\ge 1$.  \qed
\end{remark}

\begin{remark}
	One can inspect the proof to see that the monotonicity of
	$2^{-t}m(t)$ can be generalized to the condition that
	$c^t m(t)$ is strictly decreasing for some $c\in(0\,,1)$.
	\qed
\end{remark}

In order to prove the conjecture of Benjamini et al from Theorem
\ref{th:cap}, we first define $m_q(k):=[2^{k/q}]$ and interpolate
linearly to obtain a strictly increasing function [also called $m_q$]
on $\R_+$.  If $2^{-t}m(t)$ fails to be increasing near zero, then any
reasonable alteration near zero works because the conditions of
Theorem \ref{th:cap} only restrict the behavior of $m$ near
infinity. It is not hard to check that in this case,
$(\mathcal{L}g)(\lambda)\asymp\lambda^{-q}$ as $\lambda\downarrow
0$. From this it follows that $J(\rho)<\infty$ for some
$\rho\in\mathcal{P}(E)$ if and only if $\iint |x-y|^{-q} \,\rho(\d
x)\,\rho(\d y)<\infty$. This shows that the capacity criterion in
\eqref{eq:17} is the sharp necessary and sufficient condition.

We conclude by proving Theorem \ref{th:cap}

\begin{proof}[Proof of Theorem \ref{th:cap}]
  Let $\mathcal{F}:=\{\mathcal{F}_t\}_{t\ge0}$ denote the filtration
  such that each $\mathcal{F}_t$ is generated by all variables
  $X_j(r)$, where $j\ge 1$ and $r\in[0\,,t]$. If it is not already
  augmented in the usual way, then we need to augment $\mathcal{F}$ so
  that it satisfies the ``usual conditions'' of Dellacherie and Meyer
  \cite{DM}.
	
  Consider the events
  \begin{equation}
    U_n(t) := \left\{ B_k(t)=0\text{ for all $k\in\{1\,,
        \ldots,n\}$}\right\}.
  \end{equation}
  A first-passage time argument shows that for all
  $\mathcal{F}$-stopping times $\tau$,
  \begin{equation}\label{eq:U}
    \P_{1/2}\left( U_n(t+\tau) \, \big|\, \mathcal{F}_\tau
    \right) = 2^{-n} f_n(\tau),
  \end{equation}
  where for all $n\in\mathbf{N}\cup\{\infty\}$ and $\lambda\in\R$,
  \begin{equation}
    f_n(\lambda):=\prod_{k=1}^n
    \left(1-\e^{-m(k)|\lambda|}\right).
  \end{equation}
  Indeed, by the strong Markov property, it suffices to prove this for
  $\tau:=0$; and that is what we do next: Because
  $\{B_k(t)\}_{k=1}^\infty$ are conditionally independent given
  $\mathcal{F}_0$,
  \begin{equation}\label{eq:PrUProd}
    \P_{1/2}\left(U_n(t)\,\big|\, \mathcal{F}_0\right)=
    \prod_{k=1}^n\P_{1/2}\left(B_k(t)=0
      \ \big|\,\mathcal{F}_0\right).
  \end{equation}
  Let $\tau_k$ denote the first jump-time of the Markov process
  \begin{equation}
    s\mapsto \left( X_{m(k)}(s)\,,\ldots,X_{m(k+1)}(s)\right).
  \end{equation}
  The $\P_{1/2}$-law of $\tau_k$ is exponential with mean
  $1/m(k)$. Therefore, we obtain the following by splitting the
  probability according to whether or not $\tau_k>t$:
  \begin{equation}
    \P_{1/2}\left( B_k(t)=0\ \big|\, \mathcal{F}_0\right)
    =\e^{-m(k)t}\1_{\{B_k(0)=0\}}+\frac12\left[
      1-\e^{-m(k)t}\right].
  \end{equation}
  This and \eqref{eq:PrUProd} together imply \eqref{eq:U}.  Next we
  begin by recalling the argument of \cite[Theorem 4.3]{BHPS} for the
  necessity of the positive-capacity condition.
	
  Choose and fix some $\rho\in\mathcal{P}(E)$, and define
  \begin{equation}
    Z_n (\omega) :=\int\1_{U_n(t)}(\omega) \,\rho(\d t).
  \end{equation}
  It is easy to see that $\P_{1/2}(U_n(0))=2^{-n}$.  Moreover,
  stationarity and \eqref{eq:U} together imply
  \begin{equation}
    \E_{1/2}(Z_n^2) = \frac{1}{4^n}\iint
    f_n(t-s)\,\rho(\d t)\,\rho(\d s)=
    \frac{1}{4^n}\int (f_n*\rho)\,\d\rho.
  \end{equation}
  Therefore, the Paley--Zygmund lemma implies that
  \begin{equation}
    \liminf_{n\to\infty}\P_{1/2}\{Z_n>0\} \ge 
    \left[\int (f_\infty*\rho)\,\d\rho\right]^{-1}.
  \end{equation}
  From this it follows readily that
  \begin{equation}\label{eq:LB}
    \P_{1/2}\left\{ \exists t\in E:\, \sup_{k\ge 1}B_k(t)=0
    \right\}
    \ge\left[\inf_{\rho\in\mathcal{P}(E)}
      \int (f_\infty*\rho)\,\d\rho\right]^{-1}.
  \end{equation}
	
  For the converse bound we use a first-passage argument. Define
  $\tau(n)(\omega):= \inf\{t\in E:\ \omega\in U_n(t)\}$ where
  $\inf\varnothing:=\infty$; $\tau(n)$ is an $\mathcal{F}$-stopping
  time. Define $\rho_n\in\mathcal{P}(E)$ via
  $\rho(\bullet):=\rho_n(\bullet):= \P_{1/2}(
  \tau(n)\in\bullet\,|\,\tau(n)<\infty ).$ And consider the
  martingales $\{M_n(t)\}_{n=0}^\infty$, defined via
  \begin{equation}
    M_n := \E_{1/2}\left(Z_n\,\big|\,\mathcal{F}_t\right).
  \end{equation}
  Equation \eqref{eq:U} implies that $M_n(\tau(n)\wedge 1)
  \ge2^{-n}\int_{\tau(n)}^1 f_n(t-(\tau(n)\wedge 1))\,\rho_n(\d t)$.
  Because $\tau(n)\wedge 1=\tau(n)$, $\P_{1/2}$-a.s.\ on
  $\{\tau(n)<\infty\}$, it follows that
  \begin{equation}
    M_n(\tau(n)) \ge 2^{-n}\int_{\tau(n)}^1
    f_n(t-\tau(n) )\,\rho_n(\d t),
  \end{equation}
  $\P_{1/2}$-a.s.\ on $\{\tau(n)<\infty\}$. The trivial estimate
  \begin{equation}
    \E_{1/2}\left( M_n(\tau(n)\wedge 1)\right)
    \ge \E_{1/2}\left( M_n(\tau(n))\,;\tau(n)<\infty\right),
  \end{equation}
  together with the definition of $\rho_n$, impose the following:
  \begin{equation}\begin{split}
      &2^n\E_{1/2}\left( M_n(\tau(n)\wedge 1)\right)\\
      &\ge \E_{1/2}\left(\left.\int_{\tau(n)}^1
          f_n(t-\tau(n))\rho_n(\d t)\ \right|\,
        \tau(n)<\infty\right) \times \P_{1/2}\{\tau(n)<\infty\}\\
      &\ge\frac{1}{2}\int (f_n*\rho_n)\,\d\rho_n \times
      \P_{1/2}\{\tau(n)<\infty\}.
    \end{split}\end{equation}
  [The last inequality is an identity when $\rho_n$ is atomless.] By
  the optional stopping theorem, the left-most term is equal to
  $\E_{1/2}(M_n(1))=\E_{1/2}(Z_n)=2^{-n}$. Therefore,
  $\P_{1/2} \{\exists t\in E:\, \max_{1\le k\le n}B_k(t)=0 \}
  \le 2 [\int (f_n*\rho_n)\,\d\rho ]^{-1}$.
  
  We obtain an even smaller quantity if we replace $f_n$ by
  $f_\ell$, whenever $n\ge \ell$.  By Prohorov's theorem,
  there is a subsequence of $\{\rho_n\}$ and a probability measure
  $\rho$ on $E$ such that the subsequence converges weakly to
  $\rho$. It follows from the preceding that $\P_{1/2}\{\exists t\in
  E:\, \sup_{k\ge 1}B_k(t)=0\} \le 2[\int
  (f_\ell*\rho)\,\d\rho]^{-1}$.  Let $\ell\uparrow\infty$ and apply
  \eqref{eq:LB} and the monotone convergence theorem to find
  \begin{equation}\label{eq:bounds}
    \frac{1}{\inf_{\rho\in\mathcal{P}(E)}I(\rho)}\le
    \P_{1/2}\left\{\exists t\in E:\, \max_{k\ge 1}B_k(t)=0
    \right\}\le
    \frac{2}{\inf_{\rho\in\mathcal{P}(E)}I(\rho)},
  \end{equation}
  where $1/\inf\varnothing:=0$ and
  \begin{equation}
    I(\rho):=\iint\prod_{k=1}^\infty \left(
      1+\e^{-m(k)|t-s|}\right)\rho(\d s)\,\rho(\d t).
  \end{equation}
  Inequalities \eqref{eq:bounds} is valid for \emph{every} increasing
  function $m$ which satisfies $m(\mathbf{N})\subset \mathbf{N}$. Now
  we concentrate on the functions $m$ that satisfy the monotonicity
  requirements of Theorem \ref{th:cap}, and prove that $I(\rho)$ and
  $J(\rho)$ converge and diverge together [for those functions $m$].
  Our method is an adaptation of a reduction argument of \cite{BHPS},
  used to estimate Riesz-type products of the form
  $\prod_{k=1}^\infty(1+\exp(-m(k)\lambda))$ for $\lambda>0$ small.
	
  First of all, note that for integers $n\ge 1$ and real numbers
  $\lambda>0$,
  \begin{equation}\label{eq:new}\begin{split}
      \prod_{k=1}^n\left( 1+\e^{-m(k)\lambda}\right)
      &=1+\sum_{\substack{S\subseteq\{1,\ldots,n\}\\S\neq\varnothing}}
      	\e^{-\lambda\sum_{k\in S}m(k)}\\
      &\le
      	1+\sum_{\substack{S\subseteq\{1,\ldots,n\}\\S\neq\varnothing}}
      	\e^{-\lambda m(\max S)}.
    \end{split}\end{equation}
	
  For every $j=1,\ldots,n$ there are $2^{j-1}$ subsets
  $S\subset\{1\,,\ldots, n\}$ with $\max S=j$. Therefore,
  \begin{equation}\begin{split}
      \prod_{k=1}^n\left( 1+\e^{-m(k)\lambda}\right) &\le
      1+\sum_{j=1}^n 2^{j-1}\e^{-\lambda m(j)}
      \le 1+\int_1^\infty \e^{-\lambda m(\log_2 s)}\,\d s\\
      &=1+\int_{m(0)}^\infty \e^{-\lambda t}\,\d g(t) \le
      1+(\mathcal{L}g)(\lambda).
    \end{split}\end{equation}
  It follows that $I(\rho)\le 1+J(\rho)$.
	
  For the complementary bound we begin by expressing our product as
  follows: First, we write $\sum_{k\in S}m(k)$ as
  $\max(S)\sum_{k\in S}m(k)/m(\max S)$. Then,
  we write $m(t)=2^tf(2^t)$ where $f$ is increasing, and find that
  $\prod_{k=1}^n(1+\e^{-m(k)\lambda})$ is bounded below by
  \begin{equation}
		1+\sum_{\substack{S\subseteq\{1,\ldots,n\}
		S\neq\varnothing}}\exp\left(-\lambda m(\max S)\cdot
		\sum_{k\in S} 2^{k-\max S}\right).
  \end{equation}
  Since $\sum_{k\in S}2^{k-\max S}\le\sum_{j=0}^\infty 2^{-j}=2$,
  it follows that
  \begin{align}\nonumber
    \prod_{k=1}^n\left(1+e^{-m(k)\lambda}\right)
    &\ge 1+\sum_{\substack{S\subseteq\{1,\ldots,n\}\\ \nonumber
        S\neq\varnothing}}\e^{-2\lambda m(\max S)}
        =1+\sum_{j=1}^n 2^{j-1}\e^{-2\lambda m(j)}\\
    &\ge 1+\frac12\int_2^\infty\exp(-2\lambda m(\log_2 t))\,\d t.
    \label{eq:prodsum}
  \end{align}
  Because $tf(t)=m(\log_2 t)$ for $f$ increasing,
	\begin{equation}
		\prod_{k=1}^n\left(1+e^{-m(k)\lambda}\right) \ge
		1+\frac14\int_4^\infty\e^{-\lambda tf(t/2)}\,\d t
		\ge  1+\frac14\int_{m(2)}^\infty \e^{-\lambda s}\, \d g(s).
	\end{equation}
  If $D:=\sup\{|t-s|:\, s,t\in E\}$ then for all $\lambda\in(0\,,D)$,
  \begin{equation}\begin{split}
      \int_0^{m(2)} \e^{-\lambda s}\,\d g(s) &\le g(m(2)) =4,\\
      \int_{m(2)}^\infty \e^{-\lambda s}\, \d g(s)
      &\ge\int_{m(2)}^\infty \e^{-D s}\,\d g(s):=\frac 4C.
    \end{split}\end{equation}
  [We are \emph{defining} $C$ in this way.] Thus,
  $(\mathcal{L}g)(\lambda)\le (1+C)\int_{m(2)}^\infty \e^{-\lambda s}
  \,\d g(s)$ for all $\lambda\in(0\,,D)$, and therefore
  \begin{equation}
    \prod_{k=1}^n\left(1+\e^{-m(k)\lambda}\right)\ge
    1+\frac{(\mathcal{L}g)(\lambda)}{4(1+C)}\ge
    \frac{1+(\mathcal{L}g)(\lambda)}{4(1+C)}.
  \end{equation}
  This proves that $I(\rho)\ge (4+4C)^{-1}(1+J(\rho))$, whence the
  theorem.
\end{proof}

\begin{small}

  \vskip.4cm

  \noindent\textbf{Davar Khoshnevisan}\\
  \noindent Department of Mathematics, University of Utah,
  Salt Lake City, UT 84112-0090\\
  \noindent\emph{Email:} \texttt{davar@math.utah.edu}\\
  \noindent\emph{URL:} \texttt{http://www.math.utah.edu/\~{}davar}
  \\

  \noindent\textbf{David A. Levin}\\
  \noindent Department of Mathematics, University of Oregon,
  Eugene, OR 97403--1221\\
  \noindent\emph{Email:} \texttt{dlevin@uoregon.edu}\\
  \noindent\emph{URL:} \texttt{http://www.uoregon.edu/\~{}dlevin}
  \\

  \noindent\textbf{Pedro J. M\'endez-Hern\'andez}\\
  \noindent Escuela de Matem\'atica, Universidad de Costa Rica,
  San Pedro de Montes de Oca, Costa Rica\\
  \noindent\emph{Email:} 
  \texttt{pedro.mendez@ucr.ac.cr}\\
  \noindent\emph{URL:}
  \texttt{http://www2.emate.ucr.ac.cr/\~{}pmendez}
	
\end{small}

\begin{thebibliography}{99}

  \bibitem{AW:89} Arratia, R. and M. S. Waterman.  The
    Erd\H{o}s--R\'enyi strong law for pattern matching with a given
    proportion of mismatches, {\it Ann.\ Probab.}\ \textbf{17}{\it
      (3)} (1989) 1152--1169.
%
  \bibitem{AW:85} Arratia, Richard and Michael S. Waterman.  Critical
    phenomena in sequence matching, {\it Ann.\ Probab.}\
    \textbf{13}{\it (4)} (1985) 1236--1249.
%
  \bibitem{BHPS} Benjamini, Itai, Olle H{\"a}ggstr{\"o}m, Yuval Peres,
    and Jeffrey E. Steif.  Which properties of a random sequence are
    dynamically sensitive?  \textit{Ann.\ Probab.}\ \textbf{31}{\it
      (1)} (2003) 1--34.
%
  \bibitem{ChungErdos} Chung, Kai-Lai and Paul Erd\H{o}s.  On the
    lower limit of sums of independent random variables, \textit{Ann.\
      of Math. (2)} \textbf{48} (1947) 1003--1013.
%
  \bibitem{DM} Dellacherie, Claude and Paul-Andr\'e Meyer.
    \textit{Probabilities and Potential. B. Theory of Martingales},
    North-Holland Publishing Co., Amsterdam, 1982 (translated from the
    French by J. P. Wilson).
%
  \bibitem{ErdosRenyi} Erd{\H{o}}s, Paul and Alfr\'ed R{\'e}nyi.  On a
    new law of large numbers, \textit{J. Analyse Math.}\ \textbf{23}
    (1970) 103--111.
%
  \bibitem{ErdosRevesz} Erd{\H{o}}s, P. and P. R\'ev\'esz.  On the
    length of the longest head-run, in: Topics in Information Theory,
    Second Colloq., Keszthely (1975), North-Holland, Amsterdam,
    219--228. Colloq. Math. Soc. J\'anos Bolyai, Vol. 16, 1977.
%
  \bibitem{H:DPS} H{\"a}ggstr{\"o}m, Olle.  Dynamical percolation, in:
    {\it Microsurveys in Discrete Probability}, Princeton, NJ (1997),
    DIMACS Ser. Discrete Math.\ Theoret.\ Comput.\ Sci.\ \textbf{41},
    Amer.\ Math.\ Soc., Providence, RI, 59--74, 1998.
%
  \bibitem{HPS} H{\"a}ggstr{\"o}m, Olle, Yuval Peres, and Jeffrey
    E. Steif.  Dynamical percolation, {\it Ann. Inst. H. Poincar\'e
      Probab. Statist.}\ \textbf{33}{\it (4)} (1997) 497--528.
%
  \bibitem{Howroyd:01} Howroyd, J. D. Box and packing dimensions of
    projections and dimension profiles. {\it Math. Proc. Cambridge
      Phil. Soc.}  {\bf 130} (2001) 135--160.
%
  \bibitem{Kh} Khoshnevisan, Davar.  \textit{Multiparameter
      Processes}, Springer-Verlag, NY, 2002.
%
  \bibitem{KLM} Khoshnevisan, Davar, David A. Levin, and Pedro
    J. M{\'e}ndez-Hern{\'a}ndez.  Exceptional times and invariance for
    dynamical random walks, \textit{Probab. Theory Related Fields}
    \textbf{134}{\it (3)} (2006) 383--416.
%
  \bibitem{KX} Khoshnevisan, Davar and Yimin Xiao.  Packing-dimension
    profiles and fractional Brownian motion, \textit{Proc.\ Camb.\
      Phil.\ Soc.}\ {\bf 145} (2008) 205--213.
%
  \bibitem{Revesz} R{\'e}v{\'e}sz, P{\'a}l.  \textit{Random Walk in
      Random and Non-Random Environments}, second edition, World
    Scientific Publishing Co. Pte. Ltd., Hackensack, NJ, 2005.
%
\bibitem{RC} Rusakov, O. V. and A. N.  Chuprunov.
	Limit theorems for a nonhomogeneous Ornstein-Uhlenbeck process,
	\textit{Zap.\ Nauchn.\ Sem.\ S.-Peterburg. Otdel.\ Mat.\ Inst.\ Steklov} 
	\textbf{339}, Veroyatn. i Stat. 10 (2006) 111--134, 178; translation in 
	J. Math.\ Sci.\ (N.Y.) \textbf{145}{\it (2)} (2007) 4900--4913.
%
\bibitem{Rusakov} Rusakov, O. V.  A functional limit theorem for
	random variables with strong residual dependence, \textit{Teor.\
	Veroyatnost.\ i Primenen.}\ \textbf{40}{\it (4)} (1995)
	813--832.
%
\bibitem{Rusakov:94} Rusakov, O. V. Asymptotic behavior of the variance 
	of sums of random variables with random replacements (in Russian), 
	\textit{Zap.\ Nauchn.\ Sem.\ S.-Peterburg. Otdel.\ Mat.\ Inst.\ Steklov}
	\textbf{216} (1994), Problemy Teorii Veroyatnost.\ Raspred.\ \textbf{13}, 124--143, 165; 
	translation in J. Math.\ Sci.\ (N.Y.) \textbf{88}{\it (1)} (1998) 86--98.
%
\bibitem{Rusakov:93} Rusakov, O. V. 
	Central limit theorem for a summation scheme,
	in: \textit{Limit Theorems of Probability Theory}, No. 3 (in Russian), 
	241--251, 263, Izd.\ St.-Peterbg.\ Univ., St. Petersburg, 1993. 

  \end{thebibliography}
\end{document}